\documentclass[12pt]{amsart}
\usepackage{epsfig,amssymb,epic,eepic,color}
\usepackage[all]{xy}
\usepackage[T1]{fontenc}   
\newtheorem{thm}{Theorem}[section]

\newtheorem{prop}[thm]{Proposition}
\newtheorem{defn}[thm]{Definition}
\newtheorem{lemme}[thm]{Lemma}

\newtheorem{remarque}[thm]{Remark}

\newtheorem{rien}[thm]{}

\newcommand{\be}{\begin{enumerate}}
\newcommand{\ee}{\end{enumerate}}
\newcommand{\bi}{\begin{itemize}}
\newcommand{\ei}{\end{itemize}}

\def\R{\mathbb{R}}
\def\N{\mathbb{N}}

\def\H{\mathbb{H}}
\def\Z{\mathbb{Z}}

\def\E{\mathbb{E}}

\def\om{\omega}
\def\Om{\Omega}

\def\ga{\gamma}

\def\al{\alpha}
\def\be{\beta}

\def\vp{\varphi}

\def\la{\lambda}
\def\La{\Lambda}

\def\si{\sigma}
\def\Si{\Sigma}
\def\De{\Delta}

\def\ep{\varepsilon}

\def\nd{\noindent}

\def\proof{\nd {\bf Proof. }}

\begin{document}
\thispagestyle{empty}
\vskip 1cm
\begin{center}{\sc A note on the Chas-Sullivan product }
\end{center}
\title{}
\author{Fran\c cois Laudenbach}
\address{Laboratoire de math\'ematiques Jean Leray, UMR 6629 du CNRS,
Facult\'e des Sciences et Techniques,
Universit\'e de Nantes,
2, rue de la Houssini\`ere,
F-44322 Nantes cedex 3, France.}
\email{francois.laudenbach@univ-nantes.fr}

\keywords{Free loop space, intersection product, Chas-Sullivan product}
\subjclass[2000]{55N10, 57R19, 57R45}
\begin{abstract}  We give a finite dimensional approach to the Chas-Sullivan product on the free loop space of a manifold, orientable or not, compact or not.
\end{abstract}
\maketitle
\vskip 1.5cm
\section{Introduction}

Let $M$ be an $n$-dimensional manifold with empty boundary; it is not required to be either compact or orientable. Denote $LM=C^\infty(S^1, M)$ be its free 
loop space. In a famous paper \cite{chas},
 in the case when $M$ is an orientable and closed manifold,
M. Chas and D. Sullivan constructed a natural graded algebra structure on the
homology $H_*(LM;\Z)$, more precisely a product
$$H_i(LM;\Z) \otimes H_j(LM;\Z)\to  H_{i+j-n}(LM;\Z),
$$
in the same {\it intersection } spirit as the usual intersection product on $H_*(M;\Z)$.
But their ideas remained not completely accomplished.
A different approach has been considered by  R. Cohen and J. Jones 
in \cite{c-j}. According to their abstract, they describe
 "a realization of the Chas-Sullivan product in terms of a ring spectrum structure on the
 Thom spectrum of a certain virtual bundle over the loop space",
 a difficult technique
 indeed. Recently in \cite{f-t}, Y. F\'elix and J.-C. Thomas
have put {\it String Topology} into a  broad homotopy theoretical setting; they 
prove  that 
the {\it operations} in string topology are preserved by homotopy equivalence,
at least in the 1-connected case. On the contrary, in the present note we propose a finite dimensional approach, very close to the 
spirit of \cite{chas}, based on transversality arguments. We also
treat the case of a non-orientable manifold using local coefficients instead of $\Z$. \\

We do not think of $LM$ as a topological space but as a simplicial set (except in section 4). A $k$-simplex
in $LM$ is a smooth map $\Sigma : \De^k\times S^1\to M$ where $\De^k$ is the
 standard $k$-simplex. Let $0$ be the base point of $S^1=\R/\Z$. The evaluation map 
 $ev_0: LM\to M$ is simplicial: $ev_0(\Si)=\si$ where $\si: \De^k\to M$ 
 is the $k$-simplex of $M$ defined by 
 $$\si(t)=\Si(t,0).
 $$
 It is easy to form a chain complex based on the simplices of $LM$
 and a bi-complex based on {\it bi-simplices} (pairs of simplices). In order to take into 
 account the non-orientabiliy, we limit  ourselves to {\it small} simplices and bi-simplices
 (see \ref{small}). To  define an {\it intersection product} we consider {\it transverse} bi-simplexes. By Thom's transversality theorem, they generate a sub-bicomplex which
 has the same homology as $LM\times LM$. As we shall see, smallness (resp.  transversality property)
 will only refer to the 
 image of simplices (resp. bi-simplices) through the evaluation map. 
 Let us introduce the sub-bicomplex $LM\mathop{\times}\limits_{M}LM$ of 
 $LM\times LM$ made of {\it composable loops},  that is, pairs of loops having a common origin. Performing the composition yields a well defined map in homology
 $LM\mathop{\times}\limits_{M}M\to LM$
 The intersection product will 
 be not far from being defined on the chain level: the intersection of two {\it transverse}
  cycles in $LM$ will produce a "singular" smooth manifold $W$
  in $LM\mathop{\times}\limits_{M}LM$, hence in $LM$ after composing,
   which becomes a simplicial cycle once $W$ is triangulated. According to Whitehead \cite{w}, such a triangulation is unique up  to subdivision and isotopy. Therefore the product is well defined at the homology level.\\
 
 In last section it is shown that this definition of the free loop product 
 is not less efficient than the "infinite dimensional" one
 when calculations are performed in concrete geometric situations. Results due to different authors (M. Goresky, N. Hingston in \cite{g-h}, D. Chataur, J.-F. Le Borgne in \cite{chataur}) are rewritten in this setting.\\

 \nd {\sc Acknowledgments.} I had the chance to hear Dennis Sullivan lecturing 
on this subject
 on the occasion of his Doctorat Honoris Causa  at \'Ecole normale 
sup\'erieure  (Lyon) in Dec. 2001. Later on, 
 I had fruitful conversations with David Chataur 
 and Hossein Abbaspour who gave me more details. 
I feel indebted to all of them.
I am also grateful to Jean-Claude Thomas for comments on a preliminary version of this note and to  Antoine Touz\'e who prevented me from making a mistake in spectral sequences.\\

\section{Simplices and bi-simplices at the manifold level}

In this section we give a geometrical approach to the intersection product 
in the homology of $M$. The important point is that it is done in a way 
so that the intersection product lifts immediatly to $LM$.
 This requires an appeal to
Thom's transversality with constraints which is more powerful than the usual transversality, as it is usually done (see J.E. McClure \cite{clure} in the $PL$ case
or E. Castillo, R. Diaz \cite{diaz} in the smooth case).

\begin{rien}\label{small} {\rm The manifold $M$ under consideration  is equipped with an atlas 
$\mathcal A$ of charts. A {\it small} $k$-simplex is a smooth map
$$\si: \De^k \to M$$
whose image is contained in a chart from $\mathcal A$. For each simplex $\si$ some particular
chart $U(\si)\in \mathcal A$ is chosen once for all containing its image. A {\it small $k$-chain} is a 
linear combination $\xi=\sum n_i\si_i$ with coefficient $n_i\in \Z$ of finitely many small simplices $\si_i$. The {\it oriention twisted boundary } is defined by the following formula:
$$\partial \si=\sum_{i=0}^{k} \ep (-1)^i F_i\si,
$$where $F_i\si$ is the $i^{\rm th}$ face of $\si$ and $\ep$ is the sign of the Jacobian of  change of coordinates
 from $U(F_i\si)$ to $U(\si)$ calculated at any point of the image of $F_i\si$. 
The small chains with this boundary form a chain complex whose homology is
$H_*(M;\Z_{or})$, the homology with integral coefficients twisted by the orientation.}
\end{rien}

In the sequel, homotopy means smooth homotopy. Given a $k$-simplex $\si$, a homotopy $\si^t:\De^k\to M,\ t\in [0,1],$ 
with $\si^0=\si$ induces a homotopy $(F\si)^t$ of each face $F\si$. So the following definition 
makes sense.

\begin{defn} Given a chain $\xi=\sum n_i\si_i$, a {\it boundary preserving homotopy}
 of
 $\xi$ is a one 
parameter family $\xi^t, t\in [0,1]$, $\xi^t=\sum n_i\si_i^t$,
 where $t\mapsto \si_i^t$ is a homotopy of $\si_i$ into $M$ ($\si^0_i=\si_i$), 
and  such that, if $\si_i$ and $\si_j$ have a common face at time $t=0$
the corresponding simplices $\si_i^t$ and $\si_j^t$ still have a common face at any 
time.
\end{defn}
For instance, if $\xi$ is a cycle then $\xi^t$ is a cycle 
for every $t\in [0,1]$. When we consider such a homotopy of small chain, 
we shall limit  ourselves to the case when the homotopy is small, that is:
each summand $\si_i^t$ takes place in $U(\si_i)$. The following homotopy 
extension lemma is easy to prove.

\begin{lemme} Let $\tau$ be any simplex which  is a face of one of the 
summands of 
$\xi$. Any (small) homotopy of $\tau$ extends as a (small) boundary 
preserving homotopy
of $\xi$. The same statement holds for a family of homotopies with an 
extension 
  depending continuously on the parameter.\label{ext-hom}
\end{lemme} 

\begin{defn}\label{partial} A $(p,q)$-bi-simplex  of $M\times M$ is a pair
 $(u,v)$ where $u$ is a $p$-simplex and $v$ is a $q$-simplex of $M$. 
It is said to be small when both factors are small simplices of $M$.
\end{defn}
It is convenient to denote it
$u\times v$ like the map $u\times v:\De^p \times \De^q\to M\times M, 
(x,y)\mapsto (u(x),v(y))$. The (small) 
bi-simplices generate a bi-complex $C_{**}(M\times M)$
whose elements are bi-chains,
with two boundary operators twisted by the orientation, 
$$\partial_1(u\times v)=(\partial u)\times v,\  \partial_2(u\times v) = u\times \partial v,$$
and a total boundary operator
$$ D(u\times v)=\partial u\times v +(-1)^p u\times \partial v.$$

\begin{defn} A (small) bi-simplex $u\times v$ is said to be transverse when the map 
$u\times v$ is transverse to the diagonal $\De_M$ and when all its faces (they are bi-simplices) are also transverse to the diagonal.
\end{defn}
In the sequel, all transverse bi-simplices will be small,
 without saying it. The advantage is 
that $W=(u\times v)^{-1}(\De_M)$ is a proper orientable 
submanifold of 
 codimension $n$ (with corners) in $\De^p\times \De^q$. 
Moreover, if we use the charts
 $U(u)$ and $U(v)$ whose product contains the image of 
$u\times v$, $W$ receives  a canonical orientation.

Notice that, when $u\times v$ is a  transverse bi-simplex,
 a small homotopy of one of the factors keeps them transverse.
The  transverse bi-simplices generate 
a sub-bicomplex $C^{tr}_{**}(M\times M)$ of $C_{**}(M\times M)$. A bi-chain is said to be transverse when each of its bi-simplices is
transverse. It is said to be a {\it product} bi-chain if it has the form  $\xi\times \eta$
where both factors are chains in $M$.\\

\begin{lemme} Let $\xi\times \eta$ be a product bi-chain.
There exists an  arbitrarily $C^\infty$-small boundary
 preserving homotopy 
$\xi^t$ of $\xi^0=\xi$  such that $\xi^1\times \eta$ is
 transverse. Moreover, when  $\partial \xi \times \eta$ is 
transverse, the homotopy can be chosen stationary on 
$\partial \xi \times \eta$. The same statement holds for the
 second factor.\label{trans}
\end{lemme}

Considering the composition of loops we have in mind, it is very important 
to make this approximation among the bi-chains.\\

\proof   We first consider the case when $\xi\times \eta=u\times v$ is a 
bi-simplex. The conclusion follows from Thom's transversality theorem, 
not in  its elementary form but in the form known as the transversality 
with constraint (\cite{thom}). Indeed, only the first factor is moved for 
guaranteeing 
that bi-simplices remain  bi-simplices during the deformation.
 We argue as follows. Let $S$ be a small $n$-ball in the vector space
 containing the chart 
$U(u)$, so small that, for any $s\in S$, $u^s:=u+s$ is still contained 
in $U(u)$. We introduce the family $(u^s\times v)$ parametrized by $S$. 
It is transverse to 
$\De_M$ as well as its restriction to any face of $\De^p\times \De^q$. 
Therefore, according to Sard's lemma (used {\it \` a la Thom}), for almost 
every $s\in S$, $u^s\times v$ is 
a transverse bi-simplex.\\

For the general case, it is useful to observe that, in the above argument, 
the space $S$ of 
parameters can be chosen  arbitrarily small. When considering a product bi-chain 
$\xi\times \eta$, its $k$-skeleton is the collection of all the $k$-faces of 
the bi-simplices
appearing in $\xi\times \eta$. 
Arguing recursively, we may assume that all bi-simplices in 
the $k$-skeleton are transverse. Let $(u\times v)$ be a $(k+1)$-bi-simplex, 
which we 
endow with a family $S$ of translations in $U(u)$. 
According to lemma \ref{ext-hom}, the translation by any $s\in S$ extends 
as a boundary preserving homotopy 
of $\xi$ (ending at $\xi_s$) and it can be chosen small enough so that each
bi-simplex of the $k$-skeleton remains transverse. Therefore, 
for almost every  $s\in S$, $\xi_s$ is transverse along 
its $k$-skeleton and $u\times v$.  Repeating this process, 
we successively make  all  
$(k+1)$-bi-simplices transverse. 

For the relative version, we notice that, if $Fu$ is a face of $u$ and
$Fu\times v$ is transverse, then $u\times v$ is automatically transverse
 to $\De_M$ on a neighborhood of the domain of
$Fu\times v$. In that case we can moderate the translation by $s$
 so that it becomes stationary along $Fu\times v$. 
The general relative version follows easily.
\hfill$\diamond$\\

Let $[\xi]$ and $[\eta]$ be two classes of respective degrees $p$ and $q$ in 
$H_*(M;\Z_{or})$. According to lemma \ref{trans}, they can be represented by small cycles so that the bi-cycle 
$\xi\times\eta$ is transverse. This bi-cycle is unique up to
 {\it transverse homology};
more precisely we claim the following uniqueness lemma:

\begin{lemme} Let $\xi\times\eta$ and $\xi'\times\eta'$ 
be two transverse bi-cycles in the same bi-homology class ($[\xi]=[\xi']$ and 
$[\eta]=[\eta']$). There exists a transverse bi-chain $\Om$ whose total boundary
is $D\Om=\xi'\times\eta'-\xi\times\eta$.\label{uniq}
\end{lemme}

\proof There exists a 
cycle $\xi''$, homotopic  to 
$\xi$, such that $\xi''\times \eta$ and $\xi''\times \eta'$ are both transverse 
($\xi\times\eta$ is already transverse and this property is preserved
 by approximation of $\xi$; according to lemma \ref{trans}, 
such an approximation exists making $\xi''\times \eta'$ transverse). 
 If $\om$ is a $(p+1)$-chain such that $\partial \om=\xi -\xi''$, there is a 
boundary preserving homotopy $\om^t, t\in[0,1],$ of $\om$ relative to its 
boundary so that $\om^1 \times \eta$ is transverse. By the same argument there 
is a $(p+1)$-chain 
 $\om'$ with $\partial \om'=\eta-\eta'$, such that $\xi''\times \om'$ is 
transverse. 
 Eventually, there exists $\om''$ with $\partial \om''=\xi'-\xi''$ such 
that $\om''\times \eta'$ is transverse. By concatening the three transverse homologies
$\om^1\times \eta$, $\xi''\times \om'$ and $\om''\times \eta'$ we get a transverse homology joining both given bi-cycles.\hfill$\diamond$\\

\begin{rien} Intersection of cycles.\label{sign}
\end{rien}

We are now ready for defining the intersection of cycles. To begin with, we consider
 a transverse product
bi-chain $\xi\times \eta $ of degree $(p,q)$,
which reads as a sum of transverse bi-simplices
$$ \xi\times \eta = \sum n_{ij}\,u_i\times v_j.$$ 
Let $W_{ij}$ be the preimage of $\De_M$ by $u_i\times v_j$. As already noted,  
$W_{ij}$
is an oriented manifold with corners of codimension $n$. If 
 $F(u_i\times v_j)$ is a face
 of the bi-simplex, $FW_{ij}$ denotes the corresponding face of  $W_{ij}$.
 
 According to Whitehead \cite{w}, $W_{ij}$ can be smoothly triangulated  
by a $PL$-triangulation $T_{ij}$. Moreover if some faces have been already 
triangulated, by using the relative
 version of
 Whitehead's theorem, $T_{ij}$ 
 can be chosen so that the triangulated faces are subcomplexes. 
 If two bi-simplices have common faces 
 $F(u_i\times v_j)=F(u_k\times v_\ell)$, 
 then we have a canonical diffeomorphism 
 $FW_{ij}\to FW_{k\ell}$ which we think of as an identification. The triangulations of 
 these faces are chosen accordingly. 
We consider the chain of $\De_M\cong M$
 $$\xi\cdot\eta:=(-1)^{n(n-q)}\sum n_{ij} (u_i\times v_j)\vert 
(W_{ij},T_{ij})$$
 which is called the {\it intersection product}. Of course, as a chain it depends on the chosen triangulations. But, since the latter are unique up to subdividing and boundary preserving isotopy (that is, smooth isotopy of each simplex in $W_{ij}$ keeping the triangulation property), the ambiguity is not severe. The sign, which we call the {\it Dold sign}, will be  commented upon later.\\

 \begin{lemme} \label{bicycle}When $\xi\times \eta$ is a transverse bi-cycle, 
 the intersection product $\xi\cdot\eta$ is a cycle of degree $p+q-n$ 
 (with orientation twisted coefficients).
  If $\xi\times\eta$ is changed by a transverse homology (in the sense of lemma \ref{uniq}), $\xi\cdot\eta$ is changed by a homology. 
Finally, $[\xi]\cdot[\eta]$
 is well-defined in $H_{p+q-n}(M;\Z_{or})$. 
 \end{lemme}
 
 Note that a change of triangulation of the $W_{ij}$'s may be thought of as a special
  case  of a change by a transverse homology.\\
 
 \proof As $W_{ij}$ is an oriented proper submanifold, $(W_{ij},T_{ij})$ is a relative cycle in $\De^p\times \De^q$. Thus, the total 
 boundary of $(u_i\times v_j)\vert (W_{ij},T_{ij})$ is
 $$\bigl((\partial u_i)\times v_j\bigr)\vert(W_{ij},T_{ij})+(-1)^p\bigl(u_i\times (\partial v_j)\bigr)\vert(W_{ij},T_{ij}).$$ By summing on $ij$ we get the boundary of $\xi\cdot\eta$.
 As $\xi$ and $\eta$ are cycles, each hyperface in the latter  sum appears twice with opposite sign.
  The rest of the statement is easy to prove. ${} $\hfill$\diamond$\\

\begin{remarque} {\rm In his book \cite{dold} (Chap. VIII, \S 13.3),
 A. Dold explains that the 
chosen sign makes  the intersection product on homology and the (unsigned) cup-product on 
cohomology  fit together via the Poincar\'e duality. Another advantage of this sign is the following. Set $\H_*(M;\Z_{or})=H_{*+n}(M;\Z_{or})$. This regraded  homology 
endowed with the above intersection product  becomes a commutative ring in the 
graded sense.}\\
\end{remarque}

\section{Simplices and bi-simplices at  the free loop space level}

\begin{rien}{\rm We recall the evaluation map $ev_0:LM\to M$. A simplex $\Si:\De^k\times S^1\to M$ is said to be {\it small} when $\si=ev_0(\Si)$ is so. 
The $i^{th}$-face of $\Si$ is is obtained by restricting $\Si$ to 
$F_i\De^k\times S^1$. We have $F_i(ev_0(\Si))=ev_0(F_i\Si)$.
The {\it oriention twisted boundary} of $\Si$ is 
$$\partial \Si=\sum_{i=0}^{k} \ep (-1)^i F_i\Si,
$$ where $\ep$ is the sign of the Jacobian of change of  coordinates   from 
$U(ev_0(F_i\Si))$ to $U(ev_0(\Si))$. The small chains with this boundary 
form a sub-complex of $C_*(LM)$, the singular chain complex of $LM$, whose homology is $H_*(LM;\Z_{or})$. Indeed, this sub-complex is obtained from $C_*(LM)$ 
 by two operations which induce homotopy equivalences: subdivision and smoothing. The notion of boundary 
preserving homotopy is similar to the one given in \ref{small}.\label{small-loop}}
\end{rien}
\begin{lemme}The evaluation map has the lifting property 
for boundary preserving 
homotopy of chains with any initial chain. Moreover, if the lifting of the homotopy is given along some faces this partial lifting  can be extended to a global lifting.\label{lift}
\end{lemme}
\proof It is clear for $\De^k\times S^1\times[0,1]$ retracts onto $
\De^k\times S^1\times\{0\}\cup \De^k\times\{0\}\times [0,1]\cup F\times S^1\times[0,1]$, where $F$ is any
 union of faces in $\De^k$.\hfill$\diamond$\\
 
 A $(p,q)$-bi-simplex $u\times v$ in $LM\times LM $ is a map 
 $$\De^p\times \De^q\times S^1\to M\times M,$$
 $$(x,y,\theta)\mapsto \bigl(u(x,\theta),v(y,\theta)\bigr).$$
 It is said to be {\it transverse} when $\bigl(ev_0(u),ev_0(v)\bigr)$ 
 is a transverse bi-simplex
 of $M\times M$. In that case, the preimage of the diagonal $\De_M$ yields a submanifold with corners $W\subset
 \De^p\times \De^q$. For each $(x,y)\in W$, the  loops $u(x,-)$ and
$v(y,-)$ are composable, as they have common base points $u(x,0)=v(y,0)$.
 Therefore, taking a triangulation of $W$ and the Dold sign as in \ref{sign}, we get a $(p+q-n)$-chain of composable loops, which we call the {\it  loop intersection product}: 
 $$u * v:= (-1)^{n(n-q)} (u\times v)\vert W\times S^1.$$  Performing the composition (in some order, $u$ before $v$) yields a  
 $(p+q-n)$-chain in $LM$, called the {\it Chas-Sullivan product} or {\it loop product} 
 $u\mathop{\cdot}\limits_{CS}v$. Of course, it depends on the choice 
 of the triangulation of $W$. Here, we see that the product structure of $u\times v$ is very important; without it  we lose the entries of the composition. Thus, when making a bi-simplex transverse, it is crucial to do it by homotopy through bi-simplices. That is the reason why we used transversality with constraints.
 
 This loop product  extends linearly to the transverse bi-chains. When performing it 
 on a product
 bi-cycle, the result is a cycle in $LM$ whose homology class in
  $H_{p+q-n}(LM;\Z_{or})$ is well-defined. Strictly speaking the system of coefficients
  is $ev_0^*(\Z_{or})$, which we write $\Z_{or}$ for  short
 and we shall do  the same each time a loop space is in  question. More precisely, we have the following proposition.\\
  
 \begin{prop} Let $[\xi]\in H_p(LM;\Z_{or})$ and $[\eta]\in H_q(LM;\Z_{or})$.
These classes can be represented by cycles in $LM$ such that $\xi\times \eta$
 is a
 transverse bi-cycle. The class of $\xi\mathop{\cdot}\limits_{CS}\eta$ 
  is uniquely defined in $H_{p+q-n}(LM;\Z_{or})$.
\end{prop}
\proof Starting with arbitrary representatives of the given homology classes 
in $LM$, lemma \ref{trans} is applied: it produces an
 homotopy of  their images 
by the evaluation map making them a transverse bi-cycle in $M\times M$. The 
lifting homotopy property
(lemma \ref{lift}) allows one to make $\xi\times \eta$ a transverse bi-cycle 
in $LM\times LM$. For this representative, the loop product $\xi\mathop{\cdot}\limits_{CS}\eta$ is well defined. If another representative  $\xi'\times \eta'$ is used, one can prove that $\xi'\times \eta'$ and $\xi\times \eta$ are joined by a transverse homology, that is a transverse bi-chain in $LM\times LM$. It is nothing but a loop version of the  uniqueness lemma \ref{uniq} and it  can be deduced from the latter by applying the lifting homotopy property.
As a consequence the homology class of $\xi\mathop{\cdot}\limits_{CS}\eta$ is well defined.
\hfill$\diamond$\\

\begin{remarque} {\rm Of course the composition of smooth loops produces a piecewise smooth loop only. There  are two ways for correcting  this default. One consists of doing a smoothing (boundary preserving) homotopy. The other consists of using 
$LM=C^0(S^1,M)$ equipped with a mixed simplicial structure: a $k$-simplex will be
a continuous map $u:\De^k\times S^1\to M$ whose restriction to $\De^k\times {0}$ is smooth.}\\
\end{remarque}

In the sequel we use the following simplified notation: 
$\xi\cdot \eta :=\xi\mathop{\cdot}\limits_{CS}\eta$, which is defined when the bi-cycle
$\xi\times\eta$ is transverse, and $[\xi]\cdot [\eta] :=[\xi\mathop{\cdot}\limits_{CS}\eta]$ which is well defined. 
Actually, there is a 1-parameter family of compositions 
$-\cdot_s -$, $s\in [0,1]$, defined as follows. Two loops $u$ and $v$ are said to be 
{\it $s$-composable} when $u(s)=v(0)$; in that case the composed loop 
$u\cdot_s v$ is made of $u\vert [0,1-s]*v*u\vert [1-s,1]$. When $s=0$, it is the usual composition and when $s=1$ we have $u\cdot_1 v=v\cdot u$ (notice that if two loops are 0-composable, they are also 1-composable).

\begin{prop} \label{ass}At the homology level, the loop product is commutative up to sign. Precisely,  if 
$\xi$ and $\eta$ are  respectively a $p$-cycle and a $q$-cycle of $LM$, then
$[\xi]\cdot[\eta]=(-1)^{(p-n)(q-n)}[\eta]\cdot[\xi]$. The loop product
 is also associative.
\end{prop} 

If the regrading $\H_*=H_{*+n}$ is applied, then $\H_*(LM, \Z_{or})$ when endowed with the loop product  becomes a graded  commutative  and associative ring.\\

\proof  We assume that the bi-cycle $\xi\times\eta $ in $LM\times LM$ is transverse. 
We first prove: $$ [\xi]\cdot[\eta]=[\xi\cdot_1\eta].$$
Writing $\xi\times\eta=\sum n_{ij}u_i\times v_j $, one can make a boundary preserving homotopy of $\xi$ so that, for every $ij$, the following map
$$(x,y,s)\in \De^p\times\De^q\times[0,1] \mapsto  (u_i(x,s),v_j(y,0),s)\in M\times M\times
[0,1]
$$
is transverse to $\De_M\times [0,1]$. This transversality yields us a $(p+q-n+1)$-chain 
$\om $ in  $LM$ whose boundary is $\partial \om=\xi\cdot_1 \eta -\xi\cdot \eta$. Note that the loops in $\om$ have $u(x,0)$ as base points; therefore, even if the loops in
$u$ are orientation reversing no sign appears in the formula of $\partial \om$. 

Now, we  are reduced to proving:
$$[\xi\cdot_1\eta]=(-1)^{(p-n)(q-n)}[\eta]\cdot[\xi].$$
 On both sides of this equality the composition is the same. The only difference comes from the orientation of the manifold $W_{ij}$ associated to each transverse bi-simplex 
$u_i\times v_j$ appearing in $\xi\times\eta$. The permutation of both factors in this bi-simplex induces a change 
$(-1)^n$ in the co-orientation of the diagonal $\De_M$ and a change 
$(-1)^{pq}$  due to the order
of the factors $\De^p$ and $\De^q$ in the source. Moreover, the Dold sign  
$(-1)^{n(n-q)}$ is changed to $(-1)^{n(n-p)}$. Altogether, the sign rule follows. 

For associativity, if three chains $\xi, \eta, \zeta$ (of respective degrees $p,q,r$) are given in $LM$, $\xi\times\eta\times\zeta$ is easily defined to be transverse. If they are cycles and if the triple is transverse, the triple composition  $[\xi]\cdot[\eta]\cdot[\zeta]$ is well defined
in $H_{p+q+r-2n}(LM;\Z_{or})$. Moreover, in that situation the following facts are easily checked:
\begin{itemize}
\item $\xi\times \eta$ is transverse;
\item $(\xi *\eta)\times\zeta$ is transverse (where $*$ stands for the loop intersection product);
\item $(\xi\cdot\eta)\cdot \zeta$ coincides with $\xi\cdot\eta\cdot \zeta$
up to a canonical reparametrization of the circle.
\end{itemize}
The last item yields the associativity in homology once it is observed that the same 
is true for other bracketing. 

Now the only question is how to make $\xi\times\eta\times\zeta$
transverse  when the  is not. It is not sufficient to move one factor, as one factor in $M\times M\times M$ is not transverse to the small diagonal. It is necesary to move two factors, say $\xi\times\eta$,
but keeping the product structure, that is moving through product chains 
$\xi^t\times\eta^t$; it is again transversality with constraint.

\hfill$\diamond$\\

\begin{remarque}{\rm  Each free loop $\ga, \ \theta\mapsto \ga(\theta)$,
 gives rise to a 1-cycle $\bar\De(\ga)$
of loops by rotating the source: $$ \bar\De(\ga)(t)(\theta)=\ga(t+\theta),\ t\in S^1\, .$$ 
This map induces $\De:H_*(LM)\to H_{*+1}(LM)$, with twisted coefficients when 
$M$ is not orientable. Arguing  in a similar way as in the previous proposition
\ref{ass}, one could prove the theorem of Chas-Sullivan that $\H_*(LM)$, endowed with the loop product and $\De$, is a Batalin-Vilkovsky algebra.
}
\end{remarque}

\section{A multiplicative spectral sequence}

In this section, using our definition of the loop product, we discuss multiplicative properties   which have been stated and proved by Mark Goresky and Nancy Hingston
 in \cite{g-h}, \S 12 (up to the coefficient system). The setting is the one that R. Bott first considered  in his seminal paper \cite{bott1}, where he studied the standard $n$-sphere, and in the next paper \cite{bott2}. We summarize his results as follows. Let $M$ 
 be an $n$-dimensional closed Riemannian manifold whose 
 primitive geodesics are all 
simple loops with the same length (= 1, say). Denote by $\La$ the space of loops parametrized proportionally to  arc length. The class of regularity is not very important here. For instance take piecewise smooth loops;
a good class is the Sobolev class $H^1(S^1, M)$. Let $\ell^2: \La\to \R$ be the squared length. Bott proved that it is a nondegenerate function (now called a Morse-Bott function)
and he calculated the index of the critical points (that is, the closed geodesics). 
For $p\in \N$, let $\Lambda_p$ be the subspace of loops of length $\leq p$ 
and $\Si_p$ be the space of geodesics of length $p$. 
A geodesic in $\Si_p$ is just a primitive geodesic which is traversed $p$ times. As a manifold $\Si_0\cong M$ and, for $p\geq 1$, $\Si_p\cong UM$ where $UM$ stands 
for the unit tangent space to $M$.
 Let $\al_p$ be the index of the Hessian of $\ell^2$ at any 
point of $\Si_p$; obviously $\al_0=0$. Bott proved the iteration formula:
$$\al_p=p\al_1+(p-1)(n-1) .$$
Moreover he calculated (with $\Z_2$ coefficients) the spectral sequence derived from the filtration $\La_0\subset \La_1 \subset \ldots$ of $\La$. We are going to consider the same spectral sequence, up to some regrading.

\begin{defn} \label{multi}
A spectral sequence $\{E^r_{p,q}, d^r\}_{r\geq 1}$ is said to be  multiplicative when it is endowed 
with a product $E^r_{p,q}\otimes E^r_{p',q'}\to E^r_{p+p',q+q'}$ such that:
\begin{itemize}
\item[1)] $d^r$ is 
a derivation in the graded sense:
$$d^r(x\cdot y)= d^r(x)\cdot y +(-1)^{\vert x\vert}x\cdot d^r(y),$$
where $\vert.\vert$ stands for the total degree $\vert x\vert=p+q$
 when $x\in E^r_{p,q}$;
 \item[2)] the product on $E^{r+1}_{*,*}\cong H_*(E^r_{*,*},d^r)$  is induced by the one 
 of $E^r_{*,*}$ on the ideal $\ker d^r$.\\
\end{itemize}
\end{defn}

Following Chataur-Le Borgne in \cite{chataur}, we set $$\E^1_{p,q}= 
E^1_{p,q+n}=H_{p+q+n}(\La_p,
\La_{p-1}; ev_0^*(\Z_{or})).$$
The differential $d^1$  is the connecting homomorphism  in the long exact sequence of the triple $(\La_p,\La_{p-1},\La_{p-2})$. More generally, $d^r:E^r_{*,*}\to E^r_{*,*}$ is defined by the usual algorithm associated to the filtration
 $\cdots \La_{p-1}\subset \La_p\subset \La_{p+1}\subset\cdots\La$   
(see \cite{cleary},
 Chap. 2). Notice that page 0 exists in this spectral sequence, but there is no multiplicative structure on it.\\

\begin{prop} \label{spectral}The loop product endows  $\{\E^r_{*,*}\}_{r\geq 1}$ with a multiplicative structure.
\end{prop}

\proof  Here the system of coefficients $ev_0^*(\Z_{or})$ is understood. We first look at the desired  properties on page $r=1$.
When $\xi$ is a relative $i$-cycle of the pair $(\La_p,\La_{p-1})$, 
$\eta$ is a relative $j$-cycle of the pair $(\La_{p'},\La_{p'-1})$ and $\xi\times\eta$ is transverse, the loop intersection product $\xi*\eta$ is a $(i+j-n)$-chain of composable loops. By performing the composition we get a chain $c$ of loops of length $\leq p+p'$ 
since the length of the composed loop  is just the sum of the lengths of the entries. The boundary of $c$ is a cycle in $\La_{p+p'-1}$. The class of $[\xi]\cdot[\eta]$ is 
well defined in $H_{i+j-n}(\La_{p+p'},\La_{p+p'-1})$. After regrading, this product 
behaves  well with respect to   the bi-degree. It remains to check that $d^1$ is a derivation.

We look first at a transverse bi-simplex $u\times v$ of bi-degree $(i,j)$ in
 $\La\times \La$. Let $W\subset\De^i\times \De^j$ be the preimage of the diagonal
 $\De_M$ by $ev_0(u)\times ev_0(v)$. Its boundary is made of two parts:
 \begin{itemize}
 \item $\partial_1 W:=W\cap (\partial \De^i\times\De^j)$
 \item $\partial_2 W:=W\cap (\De^i\times \partial\De^j)$
 \end{itemize}
 Each part gives rise to an $(i+j+n-1)$-chain in $\La$ which is a part of 
 $\partial(u\cdot v)$. 
As in the proof of lemma \ref{bicycle} the second chain is endowed with the sign
$(-1)^i$ according to the formula for the total boundary which is written after
definition \ref{partial}. Taking the Dold sign into account yields:
$$\partial (u\cdot v)=(\partial u)\cdot v+(-1)^{n-i}u\cdot (\partial v).
$$
By summing such a formula over all bi-chains forming
$\xi \times\eta$, we get:
$$\partial (\xi\cdot \eta)=(\partial \xi)\cdot \eta+(-1)^{n-i}\xi\cdot (\partial \eta).
$$After regrading, it becomes a derivation formula. When $\xi$ and $\eta$ are relative cycles as above, $d^1([\xi]\cdot[\eta])=(d^1[\xi])\cdot[\eta]+(-1)^{\vert\xi\vert}
[\xi]\cdot d^1[\eta]$; hence, property 1 from definition \ref{multi} holds for $r=1$.\\

The product on $\E^2_{*,*}$ is defined by taking $\xi$, a relative cycle of the pair 
$(\La_p,\La_{p-2})$, and $\eta$, a relative cycle in the pair $(\La_{p'},\La_{p'-2})$, where
$\xi\times\eta$ is transverse.  Thus  the chain of composed loops is a relative cycle
 of the 
pair $(\La_{p+p'},\La_{p+p'-2})$.  It can be checked that elements in $\ker d^1$ are represesentable by such cycles. Hence property 2 from definition \ref{multi} holds for
$r=1$. \\

The same argument applies  for every $r\geq 1 $ once one  remembers
the definition of $E^r_{p,*}$ associated to the filtration of $\La$. The product on 
$\E^r_{*,*}$ is defined by taking $\xi$, a relative cycle of the pair $(\La_p,\La_{p-r})$, and
$\eta$, a relative cycle of the pair $(\La_{p'}, \la_{p'-r})$, where $\xi\times\eta$ is transverse. The boundary operator $d^r$ is induced by  $\partial$ and 
is not affected by the value  of $r$. After regrading $d^r$, like $d^1$, is a derivation and the product which is induced on its homology is the one of $E^{r+1}_{*,*}$.
\hfill$\diamond$\\

\begin{rien} The Thom isomorphism. \label{thomis}\end{rien} There is a Morse-Bott version of the famous Morse lemma. It yields a "normal" form for a Morse-Bott function near a critical manifold. We are going to apply it to the function 
$\ell^2:\Lambda\to \R$ near  the critical manifold $\Si_p$ whose  index is  $\al_p$. If one feels uncomfortable by applying this lemma in infinite dimension, one can 
take a finite dimensional approximation of $\La$ near $\Si_p$ by considering the space
of geodesics polygons (still parametrized proportionally to arc length)
 with $Np$ edges of equal length,
where $1/N$ is less than the injectivity radius. Let $E_p^-$ be the vector bundle of
rank $\al_p$ over $\Si_p$ generated by the eigenvectors of the Hessian of $\ell^2$
(with respect to some Riemannian metric on $\La$) whose eigenvalues are negative;
 it is a sub-bundle of $T\La\vert\Si_p$. With notation borrowed from \cite{g-h}, let $\Si_p^-=exp(E_p^-)$ denote the unstable manifold of $\Si_p$ with respect to the gradient of $\ell^2$; we are mainly interested
 in its germ along the critical manifold. Finally $\La_p^-$ denotes the open set of loops 
 whose length is less than $p$. We have 
 $\Si_p^-\smallsetminus \Si_p\subset\La_p^-$ (as a consequence of Taylor expansion).\\
 
The {\it Morse-Bott lemma} states:
\begin{enumerate}\item $\La_p^-$ retracts by deformation onto $\La_{p-1}$;
\item $\La_p$ retracts by deformation onto $\La_{p-1}\cup \Si_p^-$.\\
\end{enumerate}

As a consequence, the inclusions of pairs induce the following isomorphisms:
$$H_*(\La_p, \La_{p-1})\cong H_*(\La_p, \La_p^-)\cong H_*(\Si_p^-,
 \Si_p^-\smallsetminus \Si_p).$$
Here the system of coefficients is $ev_0^*(\Z_{or})$. The orientation of
 the fibre bundle 
$E_p^-$ being twisted as $ev_0^*(\Z)$ is (see \cite{g-h} Prop. 12.2),
 we are ready to apply the Thom isomorphism:
$$ h_p: H_*(\Si_p;\Z)\cong H_{*+\al_p}(\Si_p^-, \Si_p^-\smallsetminus \Si_p;ev_0^*(\Z_{or})),$$
which can also be read, when $p>0$, 
$$ h_p: H_*( UM;\Z)\cong H_{*+\al_p}(\La_p, \La_{p-1};ev_0^*(\Z_{or}) ) .$$
Notice also that 
the Gysin morphism 
$$H_*(M;\Z_{or}) \to  H_{*+n-1}(UM;\Z)
$$ 
makes 
$H_*(UM;\Z)$ a $H_*(M;\Z_{or})$-module for the intersection product. 
 In the next proposition
the coefficients  are not  written; they are meant as we just specified them.

\begin{prop}{\rm (Goresky-Hingston, \cite{g-h} Theorem 12.5).} \label{thom-p}
The Thom isomorphisms carry the 
intersection product of $H_*(UM)$ into the loop product of 
$\oplus_{p>0}H_{*+\al_p}(\La_p,\La_{p-1})$. 
 Moreover, they carry the $H_*(M)$-module structure of 
$H_*(UM)$ to the $H_*(\La_0)$-module structure of $H_*(\La_p,\La_{p-1})$.
\end{prop}


The proof below mainly 
follows the same line as \cite{g-h}. But it is based on the notion of 
loop product that we have introduced in the previous section.\\

\proof There are several steps. 

\nd A) One can factorize the intersection product $int$ of $H_*(UM)$ in the following way:
 \[\xymatrix {H_i(UM)\otimes H_j(UM)\ar[r]^-{int_0}\ar[d]_{int}&
 H_{i+j-n}(UM\mathop{\times}\limits_{M}UM)\ar[d]^{int_1}\\
H_{i+j-(2n-1)}(UM)\ar[r]^-{=}&H_{i+j-(2n-1)}(UM\mathop{\times}\limits_{UM}UM)
}\]Indeed, starting with a transverse bi-cycle $\xi\times\eta$ we may first intersect it with the fibered product 
$UM\mathop{\times}\limits_{M}UM\subset UM\times UM$  yielding an 
$(i+j-n)$-cycle $\zeta$ (this induces the morphism $int_0$), which, in turn, is transverse
to the diagonal $\De_{UM}=UM\mathop{\times}\limits_{UM}UM$ (this  induces the morphism $int_1$).\\

\nd B) Since $ev_0\vert\Si_p$ is a smooth submersion, the fibered product 
$\Si_{p,p'}:=\Si_p\mathop{\times}\limits_{M}\Si_{p'}$ is a smooth manifold. A point of it is a pair of closed geodesics of respective lengths $p$ and $p'$ which are composable as loops. In general the composed loop is not a geodesic.

Similarly, since $\Si_p^-$ is tangent to $E_p^-$, a fibre bundle over $\Si_p$,
 then
$ev_0\vert \Si_p^-$ is also a submersion  onto $M$ near $\Si_p$. Thus 
$\Si_{p,p'}^- :=\Si_p^-\mathop{\times}\limits_{M}\Si_{p'}^-$ is smooth near $
\Si_{p,p'}$. The tangent space to 
$\Si_{p,p'}^- $ along $\Si_{p,p'}$ is the restriction of $E^-_p\times E^-_{p'}$, a vector bundle 
of rank $\al_p+\al_{p'}$. So we have a Thom isomorphism:
$$h_{p,p'}:H_*(\Si_{p,p'})\cong H_{*+\al_p+\al_{p'}}(\Si_{p,p'}^-, \Si_{p,p'}^-\smallsetminus 
\Si_{p,p'}).
$$
At the chain level of small simplices, the Thom isomorphism is generated by the following Thom "extension": take a small
simplex in the base 
of a disk bundle (which hence is trivial over the considered small simplex) and cross it with the fiber. If $\xi\times\eta$ is a bi-cycle in $\Si_p\times\Si_{p'}$ transverse to $\Si_{p,p'}$, then the Thom extension $\tilde\xi\times\tilde\eta$ is transverse to
$\Si_{p,p'}^-$ and its  intersection with $\Si_{p,p'}^-$
is the Thom extension of the intersection cycle in
the base  $\Si_{p,p'}$. This proves that the Thom isomorphism carries the loop intersection product $H_i(\Si_p)\otimes H_j(\Si_{p'})\to H_{i+j-n}(\Si_{p,p'})$
to the suitable relative version of the loop intersection product 
$H_{i+\al_p}(\Si_p^-)\otimes H_{j+\al_{p'}}(\Si_{p'}^-)\to 
H_{i+j+\al_p+\al_{p'}-n}(\Si_{p,p'}^-)$. Note that, after identification, 
 the first morphism is nothing but $int_0$ from A).\\
 
 \nd C) 
 We observe that 
  $\Si_{p+p'}$ lifts (by a section of the composition map) as a submanifold of codimension $n-1$ in $\Si_{p,p'}$. 
  Indeed, any smooth geodesic  of length $p+p'$ splits uniquely into two geodesics of respective lengths $p$ and $p'$. Conversely, $(\ga,\ga')\in  \Si_{p,p'}$ 
  belongs to this lifting of $\Si_{p+p'}$
  if and only if the initial velocities 
 $\dot\ga(0)$ and
  $\dot\ga'(0)$ are positively proportional, which is a condition of codimension $n-1$. Thus, the composition is a map of pairs:
  $$comp: (\Si_{p,p'}^-, \Si_{p,p'}^-\smallsetminus \Si_{p+p'})\to (\La_{p+p'},\La_{p+p'}
  \smallsetminus \Si_{p+p'}).
  $$
  The normal bundle to  $\Si_{p+p'}$ in $\Si_{p,p'}^-$, denoted by $\nu$, is the direct sum 
 $E_p^-\oplus E_{p'}^-\oplus \nu_{p,p'}$ where $\nu_{p,p'}$ denotes the normal bundle
 to $\Si_{p+p'}$ in $\Si_{p,p'}$. The  rank of $\nu$ is $\al_p+\al_{p'}+n-1$. According to Bott's iteration formula, it equals $\al_{p+p'}$, which is the rank of 
 $E_{p+p'}^-$.
 
 If $\xi$ is a relative cycle in $(\Si_{p,p'}^-,\Si_{p,p'}^-\smallsetminus \Si_{p,p'})$ 
 transverse to $\Si_{p+p'}$, 
 its trace in the pair $(\Si_{p,p'}^-,\Si_{p,p'}^-\smallsetminus \Si_{p+p'})$ is the Thom extension of its intersection with $\Si_{p+p'}$. In other words we have the following commutative diagram:
 \[\xymatrix {H_k(\Si_{p,p'})\ar[r]^-{int_1}\ar[d]&
 H_{k-n+1}(\Si_{p+p'})\ar[d]\\
H_{k+\al_p+\al_{p'}}(\Si_{p,p'}^-,\Si_{p,p'}^-\smallsetminus\Si_{p,p'})\ar[r]&H_{k-n+1+\al_{p+p'}}(\Si_{p,p'}^-,\Si_{p,p'}^-\smallsetminus \Si_{p+p'}),
}\] 
where the vertical arrows are the respective Thom isomorphisms. The proof
of proposition \ref{thom-p} is not yet achieved because  we are still at the level of composable loops and the composition has not yet 
 been performed.
 \\
 
 \nd D) Some difficulty comes here from the fact that the composition map 
 $comp: \Si_{p,p'}^-\to \La_{p+p'}$ could be
 singular along $\Si_{p+p'}$ in the direction
  of $\nu_{p,p'}$ since $comp$ maps $\Si_{p,p'}$ into the critical level set of $\ell$
  whose value is $p+p'$.  We are going to construct:
  \begin{itemize}
  \item[a)] a linear embedding $\vp : \nu\to T\La\vert\Si_{p+p'}$, over the identity of 
  $\Si_{p+p'}$, such that $Hess(\ell^2)\circ \varphi$ is negative definite;
  \item[b)] a homotopy from $comp$ to $exp\circ \vp\circ exp^{-1}$ among the maps of pairs
  $$(\Si_{p,p'}^-,\Si_{p,p'}^-\smallsetminus \Si_{p+p'})\to (\La_{p+p'},\La_{p+p'}\smallsetminus
  \Si_{p+p'}).$$
  \end{itemize}
  If such items exist, 
  $comp_*:H_*(\Si_{p,p'}^-,\Si_{p,p'}^-\smallsetminus  \Si_{p+p'})\to 
  H_*(\La_{p+p'},\La_{p+p'}\smallsetminus\Si_{p+p'})$
 commutes with the respective Thom isomorphisms, as it is true for $\vp_*$, and the proof of proposition \ref{thom-p} is finished.
 
We choose $\vp$ to be the identity of the factor $E_p^-\oplus E_{p'}^-$. This makes sense as $\left(E_p^-\right)_\ga$ (resp. $\left(E_p'^-\right)_{\ga'}$) embeds canonically into $\left(T\La_{p+p'}\right)_{\ga*\ga'}$, thought of as the set of vector fields tangent to $M$ along 
$\ga*\ga'$: it is  a set of tangent vector fields which are vanishing
 along $\ga'$ (resp. along $\ga$); notice that a vector in $(E_p)_\ga$ corresponds to a vector field along 
$\ga$ which vanishes at $\ga(0)$ since $UM$ is homogeneous.

Let $(\ga,\ga')
\in \Si_{p+p'}\subset \Si_{p,p'}$ and $X$ be a vector in $\nu_{(\ga,\ga')}$. This vector 
indicates an infinitesimal deformation of $(\ga,\ga')$ among the pairs of composable closed smooth geodesics, deformation which separates the directions of their initial velocities. Precisely, there is a
one-parameter family $(\ga_u,\ga'_u),\ u\in [0,\ep),$ of pairs of closed geodesics such that $\ga_u(0)=\ga(0)=\ga'(0)=\ga'_u(0)$ and ${\frac{d\ }{du}\bigl(\dot\ga_u(0)
-\dot\ga'_u(0)\bigr)_{\vert u=0}=X}$ (up to a positive scalar); here we identify
$\nu_{(\ga,\ga')}$ with the fibre $UM_{\ga(0)}$. Instead of taking this deformation which leaves the length unchanged, we consider the following shortening deformation 
made of geodesic
triangles $T_u$: the first edge in $T_u$ is $\ga_u(t), t\in [0,1-\ep]$, the second edge joins geodesically
the point $\ga_u(1-\ep)$ to the nearby point $\ga'_u(\ep)$ and the third edge is 
$\ga'_u(t), t\in [\ep, 1]$ (the triangle is parametrized proportionally to arc length). We define 
$\vp(X)$ to be the infinitesimal generator of this family. By estimating $\ell(T_u)$ it is easily seen that $Hess(\ell^2)\bigl(\vp(X)\bigr)<0$. Moreover $\vp(X)$ is independent of
$E_p^-\oplus E_{p'}^-$ since, as a tangent vector field to $M$ along the composed loop $\ga*\ga'$, it does not vanish at the junction point $\ga(1)=\ga'(0)$. This proves item a). In order to obtain 
 item b) it is sufficient to make $\epsilon$ tend to 0.

 \hfill$\diamond$\\

By taking the different regradings into account, the Thom isomorphism yields
 when $p>0$:
$$\E^1_{p,q}\cong H_{p+q+n-\al_p}(\Si_p, \Z)\cong H_{p+q+n-\al_p}(UM, \Z)
\cong \H_{p+q+n-\al_p-(2n-1)}(UM, \Z),
$$
where $\H_*(UM;\Z)$ denotes the regraded intersection ring of the $(2n-1)$-dimensional
manifold $UM$. Similarly, after regrading the Gysin morphism,
$H_*(M;\Z_{or}) \to  H_{*+n-1}(UM;\Z)$, becomes a morphism of degree 0,
$\H_*(M;\Z_{or}) \to  \H_{*}(UM;\Z)$, making 
$\H_*(UM;\Z)$ a $\H_*(M;\Z_{or})$-module for the intersection product.

 The multiplicative structure on 
$\E_{*,*}^1$ can be interpreted in terms of $\H_*(UM)$.
The following statement is due to Chataur-Le Borgne
in \cite{chataur} (up to the coefficient system).

\begin{prop} ${}$ 

\nd 1) There is an isomorphism of bigraded rings:
$$\E_{*,*}^1\cong \H_*(M;\Z_{or})\oplus \H_*(UM;\Z)[T]_{\geq 1}.
$$
Here $\E_{*,*}^1$ is endowed with the bigraded  ring structure yielded by proposition 
\ref{spectral}.  The intersection rings $\H_*(M)$ and $\H_*(UM)$ have bi-degree $(0,*)$. 
The new variable $T$ has bi-degree $(1, \al_1+n-2)$ and appears only at a positive power. Regarding $\H_*(UM)$ as a $\H_*(M)$-module, the right hand side has a well-defined ring structure.

\nd 2) The differential $d^1$ on page 1 vanishes at every place.

\nd 3) The page $\infty$ inherits the same isomorphism of bigraded rings as page 1:
 $$\E_{*,*}^\infty\cong \H_*(M;\Z_{or})\oplus \H_*(UM;\Z)[T]_{\geq 1}.$$
\end{prop}

\proof 1) We read the first page of the spectral sequence,  $\E^1_{*,*}$, {\it via}  the Thom isomorphism taking proposition \ref{thom-p} into account. For instance, we have
(without writing  the coefficients):
$$\H_0(UM)=H_{2n-1}(UM)=H_{2n-1+\al_1}(UM)\cong
E^1_{1,2n-2+\al_1}=\E^1_{1,n-2+\al_{1}}$$ and $T$ is nothing but the image of the unit $\mu\in \H_0(UM)$ (that is, the fundamental class of $UM$) through the Thom isomorphism. Once the desired ring isomorphism is specified on $T$ it extends
globally using the multiplicative property of the Thom isomorphism.\\

\nd 2) As $\H_*(\La_0)$ is a direct factor in $\H_*(\La)$, $d^1: \E_{1,*}^1\to \E^1_{0,*}$
has to vanish. In particular, $d^1(T)=0$. As $d^1$ is a derivation  (prop. \ref{spectral}), it vanishes everywhere.\\

\nd 3) As a consequence, page 2 of the spectral sequence is isomorphic to page 1 as
a bi-graded ring. Therefore the differential $d^2$ vanishes for the same reason as 
$d^1$. Proceeding  recursively through
 the successive pages yields the conclusion.\hfill$\diamond$

\vfill\eject

\end{document}